# Converting of algebraic Diophantine equations to a diagonal form with the help of generalized integer orthogonal transformation, maintaining the asymptotic behavior of the number of its integer solutions

VICTOR VOLFSON

ABSTRACT. The paper presents a new generalized integer orthogonal transformation, which consists of a well-known orthogonal transform, followed by stretching the basis vectors, maintaining the asymptotic behavior of the number of integer solutions for algebraic Diophantine equation. The author shows the properties of this transformation and he receives the algorithm for finding the matrix elements of a generalized integer orthogonal transformation for algebraic Diophantine equation of the second order to diagonal form. He investigates the cases when algebraic equation can be reduced to diagonal form with the help of the generalized integer orthogonal transformation. The article includes examples illustrating the reduction of algebraic equations of the second order to the diagonal form with the help of integer generalized orthogonal transformation and of determination asymptotics behavior of integer solutions for these equations.

1. INTRODUCTION

Asymptotic estimates of the number of integer solutions for algebraic Diophantine equations in [1-4] performed mainly for the diagonal equations corresponding to the canonical equations of surfaces. Therefore it is interesting to determine the class of converting for Diophantine equation, which would result the equation to diagonal form maintaining the asymptotic behavior of the number of its integer solutions.

The paper investigates the number of integer solutions for algebraic Diophantine equations with integer coefficients. The converting must be integer and bijective to move integer solutions of algebraic Diophantine equations to integer solutions of another Diophantine equation, and vice versa.

However, not every bijection keeps the asymptotic number of integer solutions for algebraic Diophantine equations.

Keywords: algebraic Diophantine equation, diagonal form, asymptotic estimate of the number of integer solutions, generalized integer orthogonal transformation, rotation matrix.

___

Let us answer the question. What integer transforms keep the number of integer solutions of algebraic Diophantine equation in the hypercube or hypersphere? Affine conversion movement keeps the distance between the points for an algebraic integer solutions Diophantine equation, and, consequently, the number of such solutions in the hypercube or hypersphere.

On the other hand, it is known that an orthogonal transformation, which is one of the conversions of movement, can convert an algebraic surface of the second order and in some other cases of the higher orders [5] to a canonical form.

Thus, the integer orthogonal transformation, on the one hand, has the property - to keep the number of integer solutions for algebraic Diophantine equation, and, on the other hand, has the property - to convert the equation to a diagonal form. Integer matrix of orthogonal transformation is composed only of elements: 1, -1.0. This corresponds to a rotation by the angle $\pi/2$ around the axes, which significantly reduces the possibility of using this converting.

However, we are interested in integer converting which maintains not the quantity but the asymptotic behavior of the number of integer solutions for algebraic Diophantine equation in the in the hypercube or hypersphere. We'll generalize the concept of orthogonal transformation for finding such transformations.

Usually, the orthogonal transformation means an affine transformation that keeps perpendicular basis vectors and the magnitude of their lengths. We introduce the concept of "generalized orthogonal transformation," which saves only perpendicular basis vectors. A generalized orthogonal transformation also contains a deformation (stretching of basis vectors).

The resulting transformation (an orthogonal converting with a homothety), considered in [6], is a special case of a generalized orthogonal transformation, as after orthogonal transformation uses transformation strain – homothety, that performs stretching the basis vectors in the same number of times.

Now we will consider the more general case of a generalized integer orthogonal transformation, which, after orthogonal transformation performs stretching basis vectors in a different number of times. Let us show that this integer affine transformation keeps the asymptotic behavior of the number of integer solutions of an algebraic Diophantine equation.

## 2. PROPERTIES OF GENERALIZED ORTHOGONAL CONVERSIONS

Based on the properties of affine transformations points of solutions go into points of solutions, parallel lines of solutions go into parallel lines, parallel planes of solutions go into the parallel planes and $r$- dimensional parallel planes of solutions go into $r$- dimensional parallel planes solutions.

It is known that the upper limit of the number of integer solutions for algebraic Diophantine equations of degree $k$ and variables $n$ in the cube with sides $[-N, N]$ is equal to:

$$R_n^k(N) \leq k(2N+1)^{n-1} = O(N^{n-1}). \tag{2.1}$$

Thus, having in mind (2.1), we obtain the upper bound of dimension of planes for solutions for algebraic Diophantine equations $r_B = n - 1$.

We usually look for the asymptotic behavior of solutions in a square with a side $[-N, N]$ for algebraic equations of two variables. In general, when using affine transformation, the asymptotic square becomes to a parallelogram with the area equal the area of the square multiplied by the modulus of the determinant of the affine transformation - $|det(C)|$, which is not equal to 0.

Thus, the area of the parallelogram is equal to:

$$S_p = 4N^2 \cdot |det(C)|. \tag{2.2}$$

The asymptotic behavior of solutions for algebraic Diophantine equations of variables $n$ is usually studied in the hypercube with the side $[-N, N]$.

In general, using affine transformation, the asymptotic hypercube converts to inclined hyperparallelepiped with a volume equal to the volume of a hypercube multiplied by the modulus of the determinant of the affine transformation, which is not equal to 0.

Hence based on (2.2), the volume of such hyperparallelepiped is equal to:

$$V_p = (2N)^n \cdot |det(C)|. \tag{2.3}$$

It is known if the transformation is orthogonal:

$$x_1 = c_{11}x_1' + \ldots + c_{1n}x_n' + c_1, \ldots, x_n = c_{n1}x_1' + \ldots + c_{nn}x_n' + c_n, \tag{2.4}$$

then the determinant of the transformation $|det(C)|=1$

The orthogonal transformation (2.4) maintains the lengths of the segments and angles. Therefore any hypercube, when using this transformationtion, converted to same hypercube.

Let us suppose that deformation (stretching) is taken place after orthogonal transform:

$$x'_1 = k_1 x''_1, ..., x'_n = k_n x''_n,  \qquad (2.5)$$

where $k_1, ..., k_n$ are the real numbers.

The resulting affine transformation after converting (2.4), (2.5) is:

$$x_1 = c_{11} k_1 x''_1 + ... + c_{1n} k_n x''_n + c_1, ..., x_n = c_{n1} k_1 x''_1 + ... + c_{nn} k_n x''_n + c_n. \qquad (2.6)$$

We find the module of the determinant of the resulting affine transformation $C''$ having in mind (2.6). For this purpose, we use the property of determinant that multiplication of all elements of his columns (or rows) by a number is equivalent multiplying the determinant by the same number. In this case, all elements of the first column are multiplied by the value $k_1$, etc., and all elements of the column $n$ multiplied by the value $k_n$. Therefore the determinant of the orthogonal transformation (which is equal to 1) is multiplied by the module of the product $|\prod_{i=1}^{n} k_i|$:

$$|det(C'')| = \prod_{i=1}^{n} |k_i|. \qquad (2.7)$$

Having in mind (2.5), (2.7) the hypercube with the volume - $(2N)^n$ is converted to the rectangular hyperparallelepiped with the volume:

$$V_p = (2N)^n \prod_{i=1}^{n} |k_i|, \qquad (2.8)$$

when using this affine transformation.

The module of the determinant for the matrix of this transformation (when after orthogonal transformation follows the homothety with coefficient $k$) is:

$$|det(C'')| = |k|^n. \qquad (2.9)$$

Having in mind (2.8), (2.9) the hypercube with the volume $(2N)^n$ is converted to the hypercube with the volume

$$V_p = (2N|k|)^n. \tag{2.10}$$

Now we prove the following assertion.

Assertion 1

The module of the determinant of an integer affine transformation matrix $C$ is a natural number, hence the relation is performed:

$$|det(C)| \geq 1. \tag{2.11}$$

Proof

1. The determinant of integer affine transformation matrix can only be equal to an integer as its elements are integers.

2. The determinant of the matrix for non-degenerate affine transformation is not equal to 0.

Having in mind (1) and (2) the module of the determinant of matrix for an integer affine transformation is a natural number.

Corollary 1

The density of integer solutions for an algebraic Diophantine equation is multiplied by the value $|det(C)|^{-1}$ when using the integer affine transformation and therefore the density does not increase based on assertion 1.

Proof

Having in mind (2.5) and assertion 1 any hypercube converts to some hyperparallelepiped not smaller volume when using the affine integer transformation, because $|det(C)| \geq 1$.

However based on the properties of affine transformation, all integer solutions which are in the hypercube go into integer solutions in the hyperparallelepiped, therefore the solution density is multiplied by $|det(C)|^{-1}$ and it does not increase because $|det(C)| \geq 1$.

Corollary 2

The density of integer solutions for algebraic Diophantine equation does not change when you convert it using the integer affine orthogonal transformation.

The proof of this fact follows from (2.5) and using the orthogonal transformation, as the module of the determinant of the matrix for this transformation is equal to 1.

Assertion 1 can be strengthened in the particular case.

Having in mind (2.7) the value of module of the determinant for the resulting integer affine transformation, which consists of an orthogonal transformation and subsequent deformation: $x'_i = k_i x''_i, (i = 1,...,n)$ is equal to $|det(C'')| > 1$, if the value $\prod_{i=1}^{n} |k_i| > 1$.

Therefore, the density of integer solutions for the equation is decreasing using the transformation based on Corollary 1 of Assertion 1

In this case, there is increasing the distance between the points of integer solutions for the equation, parallel lines and planes of integer solutions.

Let's consider the distance between two points $M_1(x_{11},...,x_{n1}), M_2(x_{12},...,x_{n2})$:

$$|M_1 M_2| = \sqrt{(x_{12} - x_{11})^2 + ... + (x_{n2} - x_{n1})^2}.$$

Let the point $M_1$ moves to the point $M'_1$ and the point $M_2$ moves to the point $M'_2$ after the deformation conversion $x'_1 = k_1 x_1, ..., x'_n = k_n x_n$.

Then the distance between the points $M'_1$ and $M'_2$ is equal to:

$$|M'_1 M'_2| = \sqrt{k_1^2(x_{12} - x_{11})^2 + ... + k_n^2(x_{n2} - x_{n1})^2}$$

Changing the distance before and after this conversion will be characterized by the ratio of the deformation:

$$|M'_1 M'_2| / |M_1 M_2| = \sqrt{\{k_1^2(x_{12} - x_{11})^2 + ... + k_n^2(x_{n2} - x_{n1})^2 / (x_{12} - x_{11})^2 + ... + (x_{n2} - x_{n1})^2}. \quad (2.12)$$

Thus, the change of the distance (2.12), in general depends on the coordinate of points $M_1, M_2$.

Having in mind (2.12) if we using homothety ($k = k_1 = ... = k_n$):

$$|M'_1 M'_2|/|M_1 M_2| = k, \quad (2.13)$$

i.e. changing of the distance does not depends on coordinates of points $M_1, M_2$ and the distance will increase k times in this case.

If $k > 1$, then based on (2.13) after deformation the distance between the points increases, i.e. the distance between points, parallel lines and planes also increases, which corresponds to the formula (2.9).

The number of solutions for an algebraic equation is stored inside the hypersphere using an orthogonal transformation as $|det(A)| = 1$.

The density of solutions for an equation decreases in $|det(A)|$ time ($A$ is the transformation matrix) when deformation is taken place ($x''_i = k_i x'_i, (i = 1, ..., n)$) and all $|k_j| > 1$, as stretching occurs on each axis in $|k_i|$ times.

Let us explain it with the example.

Let's consider the following Thue equation:

$$x_1^2 + 2x_1 x_2 + x_2^2 - 1 = 0. \quad (2.14)$$

The matrix of the quadratic form of equation (2.14) is:

$$\begin{pmatrix} 1 & 1 \\ 1 & 1 \end{pmatrix}.$$

Equation (2.14) can be written as:

$$(x_1 + x_2)^2 - 1 = (x_1 + x_2 + 1)(x_1 + x_2 - 1) = 0. \quad (2.15)$$

Therefore, equation (2.14) has integer solutions - two parallel lines:

$$x_2 = -x_1 - 1; x_2 = -x_1 + 1.$$

The distance between the parallel lines is equal to $\sqrt{2}$.

We define all integer solutions of the equation (2.14) in the square with the sides $[-2, 2]$: $(-2, 1); (-1, 2); (-1, 0); (0, 1), (0, -1); (1, 0); (2, -1); (2, -2)$. There are 8.

We write the characteristic equation for the matrix of the quadratic form for equation (2.14): $(1-a)^2 - 1 = (2-a)a = 0$. Characteristic numbers are equal to: $a_1 = 2, a_2 = 0$.

Hence we get the transform matrix C after the orthogonal non integer coordinate conversion (rotation through an angle $\pi/4$):

$$\begin{pmatrix} \sqrt{2}/2 & -\sqrt{2}/2 \\ \sqrt{2}/2 & \sqrt{2}/2 \end{pmatrix}$$

and subsequent using homothety ($k = \sqrt{2}$), we get the resulting integer transform matrix $C_d$:

$$\begin{pmatrix} 1 & -1 \\ 1 & 1 \end{pmatrix}$$

with $|det(C_d)| = 2$.

Therefore, the resulting integer transformation corresponds to the equations:

$$x_1' = x_1 - x_2, x_2' = x_1 + x_2 \text{ or } x_1 = 0,5x_1' + 0,5x_2'; x_2 = -0,5x_1' + 0,5x_2'. \tag{2.16}$$

We substitute (2.16) in (2.15) and obtain the equation in the new coordinates:

$$(x_1 + x_2)^2 - 1 = (x_1 + x_2 + 1)(x_1 + x_2 - 1) = (x_2' + 1)(x_2' - 1) = (x_2')^2 - 1 = 0. \tag{2.17}$$

Equation (2.17) has integer solutions - two parallel lines: $x_2' = -1, x_2' = 1$. The distance between the parallel lines is equal to 2.

The distance between two parallel lines (the integer solutions of the equation (2.14)) has increased by homothety coefficient, which corresponds to (2.13).

Using integer resulting transformation, eight of these integer solutions for the equation (2.14) moved to the next integer solutions for the equation (2.17): $(-1,-1); (-1,1); (1,1); (1,-1)$ in the square with sides $[-2, 2]$. Thus, there are only four of them.

Hence their number decreased by 2 times with determinant $|det(C_d)|$, which corresponds previously said.

If we using resulting conversion with deformation, then remaining integer solutions of the equation (2.14) went beyond the scope of the square with sides $[-2, 2]$.

On the other hand there are added six integer solutions: $(-2,1); (0,1); (2,1); (-2,-1); (0,-1); (2,-1)$ of equation (2.17), which are on the lines $x'_2 = -1, x'_2 = 1$ in the square with sides $[-2, 2]$.

These solutions are obtained by converting non-integer solutions of the equation (2.14), which are on the lines $x_2 = -x_1 + 1, x_2 = -x_1 - 1$ in the square with sides $[-2, 2]$.

Therefore, the number of integer solutions for the equation (2.17) has increased from 8 to 10 compared to equation (2.14) in the square with sides $[-2, 2]$.

Thus, in general, the number of integer solutions for algebraic Diophantine equation may even increase using the conversion to diagonal form.

The question naturally arises, how the distance will change between integer solutions for non-diagonal equation when using a general resulting transformation to reduce to diagonal form, when the strain transformation is not a homothety.

The number of integer solutions for the Diophantine equation is reduced in $|det(C_d)|$ times in the hypercube ($C_d$ - resulting transformation matrix) when using this result converting, since the density of integer solutions for the Diophantine equation is also reduced in the same times.

Therefore, the distance between integer points of solutions for non-diagonal equation increases in $(|det(C_d)|)^{1/n}$ times when using a resulting transformation of the equation to diagonal form. Recall that for the resulting conversion- $|det(C_d)| > 1$.

However, if the non-diagonal equation has only a finite number of integer solutions in the hypercube with sides $[-N, N]$, then the diagonal equation will also have a finite number of integer solutions after the said conversion (in particular case, it does not have any integer solution) in the same hypercube.

Thus, there is kept the order of the number of integer solutions for the Diophantine equation - $O(1)$ in the hypercube with sides $[-N, N]$.

If the non-diagonal equation has an infinite number of integer solutions in the hypercube with sides $[-N, N]$ located on the straight line (lines), then after this transformation the diagonal equation will also have an infinite number of integer solutions in the same hypercube, located on a line (lines), i.e there is is kept the order of the number of integer solutions for the Diophantine algebraic equation - $O(N)$ in the hypercube with sides $[-N, N]$.

If the non-diagonal equation has an infinite number of integer solutions in the hypercube with sides $[-N, N]$ located on the plane (planes), then after this transformation the diagonal equation will also have an infinite number of integer solutions in the same hypercube, located on the plane (planes), i.e. there is kept the order of the number of integer solutions for the Diophantine algebraic equation - $O(N^{r-1})$ in the hypercube with sides $[-N, N]$.

If the non-diagonal equation has an infinite number of integer solutions in the hypercube with sides $[-N, N]$ located on $r$ - dimensional plane (planes), then after this transformation the diagonal equation will also have an infinite number of integer solutions in the same hypercube, located on the dimensional plane (planes), i.e, there is kept the order of the number of integer solutions for the Diophantine algebraic equation - $O(N^{r-1})$ in the hypercube with sides $[-N, N]$.

Thus, there is kept the asymptotic behavior of integer solutions for the algebraic Diophantine equation when using the resulting conversion to the diagonal form.

3. DEFINITION OF ELEMENTS FOR THE MATRIX OF A GENERALIZED ORTHOGONAL TRANSFORMATION

First we consider the case of an equation for two variables.

The matrix of positively oriented rotation by an angle can be written for this case in the form:

$$\begin{pmatrix} 1/\sqrt{(1+tg^2(a))} & -tg(a)/\sqrt{(1+tg^2(a))} \\ tg(a)/\sqrt{(1+tg^2(a))} & 1/\sqrt{(1+tg^2(a))} \end{pmatrix} \quad (3.1)$$

Note. The transformation (3.1) is valid if $a$ is not equal to $\pi/2, -\pi/2$. These cases are not of interest to us, since the rotation matrices with these angles are already integer.

We perform a homothety transformation with a coefficient $\sqrt{(1+tg^2(a))}$ after rotation transformation (3.1).

Thus, we obtain the following matrix of resulting transformation $C_d$:

$$\begin{pmatrix} 1 & -tg(a) \\ tg(a) & 1 \end{pmatrix} \qquad (3.2)$$

A transformation with a matrix $C_d$ is an integer transformation if $tg(a)$ is an integer.

It is necessary to perform the additional homothety with a coefficient $q$ to obtain an integer resultant transformation if $tg(a)$ is a rational fraction $p/q$.

It is impossible to find an integer generalized orthogonal transformation that reduces the algebraic Diophantine equation to a diagonal form if $tg(a)$ is irrational.

Now we determine in which cases $tg(a)$ is an integer or a rational fraction.

Angles of rotation $a_1, a_2$, leading the equation:

$$a_{11}x_1^2 + 2a_{12}x_1x_2 + a_{22}x_2^2 = 0 \qquad (3.3)$$

to the diagonal form is determined from the relations:

$$tg(a_1) = a_{22} - a_{11}/2a_{12} + \sqrt{(a_{22} - a_{11})^2/4a_{12}^2 + 1},$$
$$tg(a_2) = a_{22} - a_{11}/2a_{12} - \sqrt{(a_{22} - a_{11})^2/4a_{12}^2 + 1}. \qquad (3.4)$$

Having in mind (3.4) $tg(a)$ will be an integer or a rational fraction if the following relation holds:

$$1 + (a_{22} - a_{11})^2/4a_{12}^2 = p^2/q^2, \qquad (3.5)$$

where $p/q$ is a rational fraction. Naturally $a_{12}$ is not equal 0, since otherwise equation (3.3) would already have a diagonal form.

Equation (3.5) is equivalent to the Diophantine equation:

$$(a_{22} - a_{11})^2 + (2a_{12})^2 = (2a_{12}p/q)^2, \qquad (3.6)$$

which is the second-degree Fermat equation.

Поэтому уравнение (3.5) имеет следующие решения:

$$a_{11} = (u^2 - v^2)l; a_{12} = uvl; 2a_{12}p/q = (u^2 + v^2), \qquad (3.7)$$

где $u, v, l$ - натуральные числа и $u > v$.

Therefore, equation (3.5) has the following positive integer solutions:

$$a_{11} = (u^2 - v^2)l; a_{12} = uvl; 2a_{12}p/q = (u^2 + v^2), \qquad (3.7)$$

where $u, v, l$ are positive integer numbers and $u > v$.

Based on (3.7), the minimal solution of equation (3.6) is equal to: $a_{11} = 1, a_{12} = 2, a_{22} = 4$, which corresponds to equation (3.3):

$$x_1^2 + 4x_2^2 + 4x_1 x_2 = 0. \qquad (3.8)$$

Having in mind (3.4) we obtain for equation (3.5):

$$tg(a_1) = 2; tg(a_2) = -1/2. \qquad (3.9)$$

On the basis of (3.2) for the value $a_1$ we obtain an integer transformation:

$$C_{d_1} = \begin{pmatrix} 1 & -2 \\ 2 & 1 \end{pmatrix} \qquad (3.10)$$

The coefficient of homothety for this resulting transformation is:

$$\sqrt{1 + tg(a_1)^2} = \sqrt{5}. \qquad (3.11)$$

Similarly we get the non-integer conversion for the value $a_2$:

$$C_{d_2} = \begin{pmatrix} 1 & 1/2 \\ -1/2 & 1 \end{pmatrix}. \qquad (3.12)$$

We need an additional homothety transformation with the value $k = 2$ to obtain an integer conversion from (3.12).

In this case, the resulting integer transformation takes the form:

$$C_{d_3} = \begin{pmatrix} 2 & 1 \\ -1 & 2 \end{pmatrix}. \qquad (3.13)$$

The coefficient of homothety for the resulting integer transformation (3.13) is defined as follows:

$$k\sqrt{1 + tg(a_2)^2} = 2\sqrt{1 + (1/2)^2} = \sqrt{5}, \qquad (3.14)$$

which corresponds to (3.11).

We have considered an example when $tg(a)$ is an integer, or a rational fraction.

If we use homothety, then all elements of the transformation matrix are multiplied by the homothety coefficient $k$, so the determinant of the transformation is multiplied by $k^n$, where $n$ is the number of variables.

If the resultant transformation includes a homothety with the coefficient $k$, then taking in mind that the determinant of the orthogonal transformation is equal to 1, the modulus of the determinant of this resulting transformation is equal to $|det(C_d)| = k^n$ and the coefficient of homothety is equal to:

$$k = (|det(C_d)|)^{1/n}. \tag{3.15}$$

Based on (3.15) the value of the homothety coefficient in the above example is equal to $k = (|det(C_d)|)^{1/2} = \sqrt{5}$, which corresponds to (3.11) and (3.14).

Now we consider an example when $tg(a)$ is irrational. Such cases more often found in the practice of Diophantine equations.

As mentioned before, it is impossible to find a generalized orthogonal transformation in these cases, which would lead the equation to diagonal form.

I will give a simple example of such Diophantine equation:

$$x_1^2 + x_1 x_2 = 0. \tag{3.16}$$

Having in mind (3.4) we define the following values for equation (3.16):

$$tg(a_1) = -1 - \sqrt{2}, tg(a_2) = -1 + \sqrt{2}. \tag{3.17}$$

Based on (3.17) we obtain two non-integer orthogonal transformations:

$$C_1 = \begin{pmatrix} 1 & 1+\sqrt{2} \\ -1-\sqrt{2} & 1 \end{pmatrix} \text{ and } C_2 = \begin{pmatrix} 1 & 1-\sqrt{2} \\ -1+\sqrt{2} & 1 \end{pmatrix}, \tag{3.18}$$

which can not be convert by stretching into an integer form.

Now we consider a more general case of determining the elements for the matrix of a generalized orthogonal transformation.

Assertion 2

Let the matrix of the quadratic form of an algebraic Diophantine second order equation for $n$ variables have eigenvectors with integer coordinates, respectively:

$$a_1 = (a_{11},...,a_{1n}),...,a_n = (a_{n1},...,a_{nn}), \qquad (3.19)$$

then the integer generalized orthogonal transformation that reduces the above algebraic Diophantine equation to a diagonal form has the form:

$$C = \begin{pmatrix} a_{11} & ... & a_{n1} \\ ... & & ... \\ a_{1n} & ... & a_{nn} \end{pmatrix}. \qquad (3.20)$$

Доказательство

Proof

If the eigenvectors of the matrix of the quadratic form have coordinates (3.19), then the lengths of the corresponding eigenvectors will be:

$$|a_1| = \sqrt{a_{11}^2 + ... + a_{1n}^2},...,|a_n| = \sqrt{a_{n1}^2 + ... + a_{nn}^2}. \qquad (3.21)$$

Having in mind (3.21), the vectors of the new orthonormal basis will have coordinates:

$$a'_1 = (a_{11}/\sqrt{a_{11}^2 + ... + a_{1n}^2},...,a_{1n}/\sqrt{a_{1n}^2 + ... + a_{1n}^2}),$$

...,

$$a'_n = (a_{n1}/\sqrt{a_{11}^2 + ... + a_{1n}^2},...,a_{nn}/\sqrt{a_{1n}^2 + ... + a_{1n}^2}). \qquad (3.22)$$

Based on (3.22), the matrix of the orthogonal transformation will have the form:

$$S = \begin{pmatrix} a_{11}/\sqrt{a^2_{11} + ... + a^2_{1n}} & ... & a_{1n}/\sqrt{a^2_{11} + ... + a^2_{1n}} \\ ... & ... & ... \\ a_{n1}/\sqrt{a^2_{11} + ... + a^2_{1n}} & ... & a_{nn}/\sqrt{a^2_{11} + ... + a^2_{1n}} \end{pmatrix}. \qquad (3.23)$$

The orthogonal transformation (3.23) leads the original algebraic equation to a diagonal form, but in the general case it is not an integer transformation, since in the denominator of the elements of matrix there may be either irrational numbers or integers that do not divide the numerator without remnant.

We perform the transformation of the deformation additionally, which maintains the orthogonality of the basis vectors:

$$x'_1 = k_1 x''_1, ..., x'_n = k_n x''_n,  \qquad (3.24)$$

where $k_1 = \sqrt{a_{11}^2 + ... + a_{1n}^2}, ..., k_n = \sqrt{a_{n1}^2 + ... + a_{nn}^2}$.

After the additional transformation (3.24) we obtain the generalized orthogonal transformation:

$$C = \begin{pmatrix} a_{11} & ... & a_{1n} \\ ... & & ... \\ a_{n1} & ... & a_{nn} \end{pmatrix},$$

which has integer coefficients, based on the conditions of the assertion.

Let us explain assertion 2.

It is known that, in the general case, an algebraic Diophantine equation of second order for $n$ variables with integer coefficients:

$$F(x_1, ..., x_n) = \sum_{i,j=1}^{n} a_{ij} x_i x_j + 2 \sum_{i=1}^{n} x_i + a_0 = 0 \qquad (3.25)$$

using the transformation of the origin:

$$x_1 = x'_1 + x_{10}, ..., x_n = x'_n + x_{n0}. \qquad (3.26)$$

can be represented as:

- for the central case:

$$F(x'_1, ..., x'_n) = \sum_{i,j=1}^{n} a_{ij} x'_i x'_j + a'_0 = 0, \qquad (3.27)$$

- for the non-center case:

$$F(x'_1, ..., x'_n) = \sum_{i,j=1}^{n-1} a_{ij} x'_i x'_j + 2 a'_n x'_n = 0. \qquad (3.28)$$

where $a'_0 = F(x_{10}, ... x_{n0})$ $a'_0 = F(x_{10}, ... x_{n0})$.

Since the function $F(x_1, ..., x_n)$ in equation (3.25) is integer, the coefficients $a'_0, a'_k$ are also integer and the coefficients of the quadratic form $a_{ij}$ of equation (3.25) do not change at all if the transformation of the origin (3.26) is also integer.

We assume that there is an integer transformation of the origin, which converts an algebraic Diophantine equation of the second order to equations (3.27) or (3.28) with integer coefficients in the statement 2.

The corollary of Assertion 2

Let the matrix of the quadratic form $F(x_1,...,x_n)$ have eigenvalues: $t_1,...,t_n$ and eigenvectors with integer coordinates: $(c_{11},...,c_{1n}),...,(c_{n1},...,c_{nn})$.

Then a homogeneous Diophantine equation of the second order $F(x_1,...x_n)=0$ is reduced to a diagonal homogeneous Diophantine equation with integer coefficients:

$$F(x'_1,...,x'_n) = (\sum_{i=1}^{n} c_{1i}^2)t_1(x'_1)^2 + ... + (\sum_{i=1}^{n} c_{ni}^2)t_n(x'_n)^2 = 0 \qquad (3.29)$$

using the generalized orthogonal transformation with a matrix:

$$C = \begin{pmatrix} c_{11} & ... & c_{n1} \\ ... & & ... \\ c_{1n} & ... & c_{nn} \end{pmatrix}.$$

Proof

Having in mind Assertion 2 using the generalized orthogonal transformation:

$$C = \begin{pmatrix} c_{11} & ... & c_{n1} \\ ... & & ... \\ c_{1n} & ... & c_{nn} \end{pmatrix}$$

a homogeneous Diophantine equation of the second order $F(x_1,...x_n)=0$ is reduced to a homogeneous diagonal Diophantine equation with integer coefficients:

$$F(x'_1,...,x'_n) = k_1^2 t_1 (x'_1)^2 + ... + k_n^2 t_n (x'_n)^2, \qquad (3.30)$$

where $t_1,...,t_n$ are the eigenvalues of the matrix for the quadratic form $F(x_1,...x_n)$, and $k_j (1 \leq j \leq n)$ is the deformation coefficient by the $j$- coordinate which is determined as:

$$k_j = \sqrt{\sum_{i=1}^{n} c_{ji}^2}. \qquad (3.31)$$

Подставим формулу (3.31) в (3.30) и получим (3.29).

We substitute the formula (3.31) into (3.30) and get (3.29).

Now we'll consider at an illustrative example.

It is necessary to determine the number of integer solutions for the Diophantine equation:

$$(x_1)^2 + 5(x_2)^2 + (x_3)^2 + 2x_1x_2 + 6x_1x_3 + 2x_2x_3 + 8x_1 + 20x_2 + 16 = 0. \tag{3.32}$$

The matrix of the quadratic form of equation (3.32) has the form:

$$A = \begin{pmatrix} 1 & 2 & 3 \\ 1 & 5 & 1 \\ 3 & 1 & 1 \end{pmatrix},$$

where $det(A) = -36$ (is not equal to 0).

Consequently, equation (3.32) corresponds to the central surface. Let us find the coordinates of its center from the system of equations:

$$x_{10} + x_{20} + 3x_{30} = 4, \; x_{10} + 5x_{20} + x_{30} = 10, \; 3x_{10} + x_{20} + x_{30} = 0. \tag{3.33}$$

The system (3.33) has a unique solution: $x_{10} = -1, x_{20} = 2, x_{30} = 1$ since the determinant of system is not equal to 0.

Consequently, it is possible an integer transfer of the origin to the center of the surface, therefore assertion 2 holds for the equation (3.32).

We introduce new coordinates and transfer the center of the surface to the origin of the new coordinates of the system:

$$x_1 = x_1' + 1, x_2 = x_2' - 2, x_3 = x_3' - 1. \tag{3.34}$$

We substitute the new coordinates (3.34) in (3.32) and obtain the equation:

$$(x_1'+1)^2 + 5(x_2'-2)^2 + (x_3'-1)^2 + 2(x_1'+1)(x_2'-2) + \\ 6(x_1'+1)(x_3'-1) + 2(x_2'-2)(x_3'-1) + 8(x_1'+1) + 20(x_2'-2) + 16. \tag{3.35}$$

We obtain the equation in new coordinates by giving similar terms in (3.35):

$$(x_1')^2 + 5(x_2')^2 + (x_3')^2 + 2x_1'x_2' + 6x_1'x_3' + 2x_2'x_3' = 0. \tag{3.36}$$

Let us find the eigenvalues of the matrix for the quadratic form of equation (3.36):

$$\begin{bmatrix} 1-t & 2 & 3 \\ 1 & 5-t & 1 \\ 3 & 1 & 1-t \end{bmatrix} = -t^3 + 7t^2 - 36 = 0. \tag{3.37}$$

Based on the characteristic equation (3.37) we obtain the following eigenvalues: $t_1 = 3, t_2 = 6, t_3 = -2$.

Now we find an eigenvector for the eigenvalue $t_1 = 3$, solving the system of equations:

$$-2a_1 + a_2 + 3a_3 = 0, a_1 + 2a_2 + a_3 = 0, 3a_1 + a_2 - 2a_3 = 0. \tag{3.38}$$

The system of linear equations (3.38) has a solution: $a_1 = 1, a_2 = -1, a_3 = 1$. Therefore, we take the integer vector:

$$a = \begin{pmatrix} 1 \\ -1 \\ 1 \end{pmatrix}$$

as an eigenvector.

Let us find the eigenvector for the eigenvalue $t_2 = 6$, solving the system of equations:

$$-5b_1 + b_2 + 3b_3 = 0, b_1 - b_2 + b_3 = 0, 3b_1 + b_2 - 5b_3 = 0. \tag{3.39}$$

The system of linear equations (3.39) has a solution: $b_1 = 1, b_2 = 2, b_3 = 1$, Therefore, we take the integer vector:

$$b = \begin{pmatrix} 1 \\ 2 \\ 1 \end{pmatrix}$$

as the eigenvector.

Let us find the eigenvector for the eigenvalue $t_3 = -2$, solving the system of equations:

$$3c_{11} + c_2 + 3c_3 = 0, b_1 + 7c_2 + c_3 = 0, 3c_1 + c_2 - 3c_3 = 0. \tag{3.40}$$

The system of linear equations (3.40) has a solution: $b_1 = -1, b_2 = 0, b_3 = 1$. Therefore, we take the integer vector:

$$c = \begin{pmatrix} -1 \\ 0 \\ 1 \end{pmatrix}$$

as the eigenvector.

Based on Assertion 2 the integer generalized orthogonal transformation leading (3.36) to the diagonal form for this example is:

$$C = \begin{pmatrix} 1 & 1 & -1 \\ -1 & 2 & 0 \\ 1 & 1 & 1 \end{pmatrix}. \tag{3.41}$$

Having in mind (3.41), the value $\det(C) = 6$ (not equal to 0), therefore the transformation $C$ is nondegenerate.

Based on (3.32), the diagonal Diophantine equation (obtained after the generalized orthogonal transformation) will have the form:

$$F(x'_1,...,x'_n) = 9(x'_1)^2 + 36(x'_2)^2 - 4(x'_3)^2 = (3x'_1)^2 + (6x'_2)^2 - (2x'_3)^2 = 0. \tag{3.42}$$

Equation (3.42) is the second-order Fermat equation for the variables: $3x'_1, 6x'_2, 2x'_3$. The positive integer solutions of equation (3.42) can be written in the form:

$$3x' = (m^2 - n^2)l, \ 6x'_2 = 2mnl, \ 2x'_3 = (m^2 + n^2)l, \tag{3.43}$$

where $m, n, l \, (m > n)$ are natural numbers.

Having in mind (3.43) a maximum value has the variable:

$$2x'_3 = (m^2 + n^2)l \text{ or } x'_3 = (m^2 + n^2)l/2. \tag{3.44}$$

The value $(m^2 + n^2)l/2$ is a natural number, if $(m^2 + n^2)l$ is even. The latter is possible in cases when either $l$ - even, or $m^2 + n^2$ - even. The value $m^2 + n^2$ will be even if $m$ and $n$ are even, or $m$ and $n$ are odd. It will be in half cases of values $m, n$.

Thus, $(m^2 + n^2)l$ will be even in $3/4$ cases of values $m, n, l$.

Then, based on (3.44), the number of natural solutions of equation (3.42) in a cube with side $N$ is found from the condition:

$$(m^2 + n^2)l \leq 2N. \tag{3.45}$$

Therefore, we must find the number of points with natural coordinates satisfying condition (3.45), for which $(m^2 + n^2)l$ is even.

Having in mind (3.45) if $l = 1$, we obtain the inequality $m^2 + n^2 \leq 2N$. Since $m > n$, then we get that it is necessary to find the number of points with positive integer coordinates inside a sector with a radius equal to $(2N)^{1/2}$ above the main diagonal.

Based on [7] the number of such points will be:

$$\pi N / 4 + O(N^{1/2}). \tag{3.46}$$

Having in mind (3.45), if $l = l_{max}$, we obtain the inequality:

$$m^2 + n^2 \leq 2N / l_{max}. \tag{3.47}$$

The number of points with positive integer coordinates satisfying condition (3.47) is obtained similarly:

$$\pi N / 4l_{max} + O(N^{1/2}). \tag{3.48}$$

Based on (3.43), (3.45) we obtain the following upper bound for all cases:

$$\pi N / 4(1 + 1/2 + \ldots + 1/l_{max}) + O(N^{1/2}). \tag{3.49}$$

The function $1/k$ is strictly decreasing and, based on the Euler-McLeron formula, we obtain:

$$\sum_{k=1}^{l_{max}} 1/k = \int_{k=1}^{l_{max}} dx / x + C + O(1/l_{max}), \tag{3.50}$$

where $C$ is a constant.

We substitute (3.50) into (3.49) and obtain the following upper bound:

$$\pi N / 4[\ln(l_{max}) + C + O(1/l_{max})] + O(N^{1/2}). \tag{3.51}$$

Since $l_{max} < N$, based on (3.51), and taking into account $3/4$ cases for which $(m^2 + n^2)l$ are even, we obtain an upper estimate for the number of natural solutions for equation (3.42):

$$R_3^+(N) \leq 3\pi N \ln(N) / 16 + O(N). \tag{3.52}$$

The function $F(x_1', x_2', x_3')$ in equation (3.42) is even for all the variables, so the integer solution of equation (3.42) may be in any of the eight octants of a cube with sides $[-N, N]$ compared to one octant for natural values.

Therefore, we must multiply (3.52) by the value 8 to obtain an estimate of the number of integer solutions for equation (3.42):

$$R_3(N) \leq 3\pi N \ln(N)/2 + O(N). \tag{3.53}$$

Having in mind, that $3\pi N \ln(N)/2 = O(N\ln(N))$, the expression (3.53) can be written in the form:

$$R_3(N) \ll N\ln(N). \tag{3.54}$$

Since the asymptotic estimate of the number of integer solutions for an algebraic Diophantine equation are saved when using the integral generalized orthogonal transformation, our estimate (3.54) is true for the number of integer solutions for the original Diophantine equation of our example (3.32).

The Diophantine equation (3.42) is homogeneous, so the point $(0,0)$ is its solution. In addition, there are integer solutions (3.42) on two straight lines that are in the plane $x'_1 = 0$:

$$3x'_1 + x'_3 = 0, 3x'_1 - x'_3 = 0, \tag{3.55}$$

as well as the integer points of two straight lines in the plane $x'_2 = 0$:

$$3x'_1 + 2x'_3 = 0, 3x'_1 - 2x'_3 = 0. \tag{3.56}$$

The asymptotic estimate of the integer solutions (3.55) and (3.56) does not exceed $O(N)$, therefore the asymptotic estimate for the number of integer solutions for Diophantine equation (3.42) and equation (3.32) remain unchanged and is determined by (3.53) or (3.54) respectively.

All integer solutions for equations (3.42) and (3.32) are on rectilinear generators of the corresponding cones. Based on estimate (3.54), the upper bound for the number of such generators is $O(\ln(N))$.

The transformation of algebraic Diophantine equations to a diagonal form is also a method of solving it, since finding the integer solutions of equation (3.32) is considerably more complicated than for equation (3.42).

Let's consider the following two points of the natural solutions of equation (3.42):

1) values $l=6, m=2, n=1$ corresponding to the point $x'_1 = 6, x'_2 = 4, x'_3 = 15$ and integer points on the rectilinear generator $x'_1 = 6t_1, x'_2 = 4t_1, x'_3 = 15t_1$, where the variable $t_1$ takes any integer values;

2) values $l=6, m=3, n=1$ corresponding to the point $x'_1 = 16, x'_2 = 6, x'_3 = 30$ and integer points on the rectilinear generator $x'_1 = 16t_2, x'_2 = 6t_2, x'_3 = 30t_2$, where the variable takes any integer values.

It is easy to find the corresponding solutions of equation (3.32) using the transformation:

$$x_1 = x'_1 + x'_2 - x'_3 + 1, x_2 = -x'_1 + 2x'_2 - 2, x_3 = x'_1 + x'_2 + x'_3 - 1 \qquad (3.57)$$

For example, based on (3.57) the point 1) of the integer solutions for equation (3.32) corresponds to the following solution point of equation (3.42): $x_1 = -4, x_2 = 0, x_3 = 24$, thus rectilinear generator of the cone: $x_1 = -4t_3, x_2 = 0, x_3 = 24t_3$, and point 2) of the integer solutions for equation (3.32) corresponds to the following solution point of equation (3.42): $x_1 = -7, x_2 = -6, x_3 = 51$, thus rectilinear generator of the cone: $x_1 = -7t_4, x_2 = -6t_4, x_3 = 51t_4$.

One can obtain other integer solutions of the Diophantine equation (3.32) in a similar way.

## 4. CONCLUSION AND SUGGESTIONS FOR FURTHER WORK

The next article will study integer non-orthogonal transformations that maintain the asymptotic behavior of the number of integer solutions for the algebraic Diophantine equation.

## 5. ACKNOWLEDGEMENTS

Thanks to everyone who has contributed to the discussion of this paper.